%% file: ADL-draft.tex
\documentclass[12pt]{elsarticle}
\usepackage[T1]{fontenc}
\usepackage{amsmath,amssymb,xcolor,url,bbm,mathtools}
\usepackage{graphicx}
\usepackage{float}
\usepackage{subcaption}
\usepackage[normalem]{ulem}
\usepackage{algorithm}
\usepackage{algpseudocode}

\input{notation}

\begin{document}

\begin{frontmatter}
\title{Approximate Deconvolution Leray Reduced Order Model for Convection-Dominated Flows}
\author[trento]{Anna Sanfilippo}
\author[vt]{Ian Moore}
\author[cattolica]{Francesco Ballarin}
\author[vt]{Traian Iliescu}

\affiliation[trento]{
organization={Department of Mathematics, Universit{\`a} di Trento},
addressline={via Sommarive 14}, 
city={Povo},
postcode={38123}, 
country={Italy}
}
\affiliation[vt]{
organization={Department of Mathematics, Virginia Tech},
addressline={225 Stanger Street}, 
city={Blacksburg},
postcode={24061}, 
state={Virginia},
country={US}
}
\affiliation[cattolica]{
organization={Department of Mathematics and Physics, Universit{\`a} Cattolica del Sacro Cuore},
addressline={via Garzetta 48}, 
city={Brescia},
postcode={25133}, 
country={Italy}
}

\begin{abstract}
In this paper, we propose a novel ROM stabilization strategy for under-resolved convection-dominated flows, the approximate deconvolution Leray ROM (ADL-ROM).
The new ADL-ROM introduces AD as a new means to increase the accuracy of the classical Leray ROM (L-ROM) without degrading its numerical stability.
We also introduce two new AD ROM strategies: the Tikhonov and van Cittert methods.
Our numerical investigation for convection-dominated systems shows that, when the filter radius is relatively large, the new ADL-ROM is more accurate than the standard L-ROM. 
Furthermore, the new ADL-ROM is less sensitive with respect to model parameters than L-ROM.
\end{abstract}

\begin{keyword}
Reduced order models \sep approximate deconvolution \sep under-resolved regime \sep spatial filter \sep regularization \sep Leray model



\end{keyword}
\end{frontmatter}

\section{Introduction}

Galerkin ROMs (G-ROMs) are computational models that leverage data to dramatically reduce the dimension of full order models (FOMs), i.e., models obtained from classical numerical discretizations (e.g., the finite element method (FEM)).
G-ROMs have been used to reduce the FOM computational cost 
in the numerical simulation of laminar fluid flows described by the Navier-Stokes equations (NSE)~\cite{ballarin2015supremizer,brunton2019data,hesthaven2015certified,noack2005need,quarteroni2015reduced}.
However, in the under-resolved regime (i.e., when the number of ROM degrees of freedom is not enough to accurately represent the flow dynamics), G-ROM yields inaccurate results, usually in the form of numerical oscillations.
To alleviate this numerical inaccuracy, several types of ROM  stabilization techniques have been developed~\cite{AliBallarinRozza2020,AliBallarinRozza2021,balajewicz2012stabilization,barone2009stable,bergmann2009enablers,carlberg2011efficient,fick2018stabilized,grimberg2020stability,huang2022model,kalashnikova2010stability,pacciarini2014stabilized,parish2020adjoint,parish2023residual,torlo2017stabilized,zoccolan2023streamline}.

One popular type of ROM stabilization for fluid flows 
is {\it regularized ROMs (Reg-ROMs)} 
\cite{Girfoglio2019CF,girfoglio2021pod,gunzburger2020leray,iliescu2018regularized,kaneko2020towards,rezaian2023predictive,sabetghadam2012alpha,tsai2022parametric,wells2017evolve,xie2018numerical}. 
The principle used to construct Reg-ROMs is to use a {\it ROM  filter} to smooth out various terms in the underlying NSE.
There are two ROM filters in current use:
(i) The ROM projection is a ROM filter in which filtering takes place exclusively in the ROM space.
(ii) The ROM differential filter is a ROM filter in which filtering takes place exclusively in the physical space.
Higher-order DFs have been introduced in~\cite{gunzburger2019evolve}.
Both the ROM projection and the ROM differential filter have been used to construct Reg-ROMs.

One of the most popular Reg-ROMs is the {\it Leray ROM (L-ROM)}~\cite{Girfoglio2019CF,girfoglio2021pod,kaneko2020towards,sabetghadam2012alpha,wells2017evolve}, which is a ROM stabilization for under-resolved convection-dominated flows.
The L-ROM, inspired by the work of Jean Leray in the mathematical theory of the NSE~\cite{leray1934sur}, is based on a simple, yet powerful idea:
Replace the nonlinear term in the NSE, $\bu \cdot \nabla \bu$, with $\obu \cdot \nabla \bu$, where $\obu$ denotes the ROM filtered velocity.
Because the filtered velocity, $\obu$, is smoother than the original, unfiltered velocity, $\bu$, Leray was able to prove the existence of weak solutions to the NSE~\cite{leray1934sur}.
Filtering the velocity is also a good idea from the computational point of view:
The Leray model was used in under-resolved  simulations of turbulent flows with classical numerical discretizations~\cite{geurts2003regularization}.
The Leray model was extended to reduced order modeling of the NSE in~\cite{wells2017evolve} (see~\cite{sabetghadam2012alpha} for earlier work on the Kuramoto-Sivashinsky equations). 
In~\cite{wells2017evolve}, it was shown that the Leray ROM yields significantly more accurate results than the standard Galerkin ROM in the under-resolved numerical simulation of the 3D flow past a circular cylinder at Reynolds number $Re=1000$.
Since then, the L-ROM has been successfully used in the under-resolved simulation of various convection-dominated flows, e.g., further NSE applications~\cite{Girfoglio2019CF,girfoglio2021pod,girfoglio2021pressure,Girfoglio2023JCP}, the stochastic NSE~\cite{gunzburger2019evolve,gunzburger2017ensemble}, and the quasi-geostrophic equations~\cite{Girfoglio2023CRM,Girfoglio2023JCAM}.

Despite the significant improvement over the standard Galerkin ROM, one drawback of the L-ROM observed in these numerical investigations is that it can be overdiffusive.
For example, when the ROM filter radius, $\delta$, is too large, the filter introduces too much dissipation in the L-ROM.
Another drawback of the L-ROM observed in these numerical investigations is its sensitivity with respect to the ROM filter radius.
Specifically, small perturbations of the optimal filter radius can decrease the L-ROM accuracy.

To address these L-ROM drawbacks, we propose a new Reg-ROM, the {\it approximate deconvolution Leray ROM (ADL-ROM)}, which is a significant improvement of the L-ROM.
The innovation in the new ADL-ROM is the leverage of approximate deconvolution, a classical strategy used in the image processing and inverse problems communities, to increase the accuracy and decrease the sensitivity of the L-ROM. 
Specifically, in the new ADL-ROM, the filtered L-ROM velocity, $\obu$, is replaced with the AD velocity, $D(\obu)$. 
Thus, the ADL-ROM nonlinear term is $D(\obu) \cdot \nabla \bu$.
Since $D(\obu)$ approximation of $\bu$ is both stable and accurate, using the AD operator in the ADL-ROM has the role of increasing the accuracy of the nonlinear term while preserving the extra stability introduced by the ROM filtering in $\obu$.
Thus, the ADL-ROM can be thought of as a compromise between the  G-ROM and the overdissipative L-ROM.
To our knowledge, this is the first time the AD concept is used to construct ROM stabilizations, such as Reg-ROMs.

The rest of the paper is organized as follows:
In Section~\ref{sec:g-rom}, we outline the construction of the standard G-ROM for the NSE.
In Section~\ref{sec:ad-rom}, we first present the differential filter and then introduce several  AD strategies: the van Cittert, Tikhonov, and Lavrentiev approaches.
In Section~\ref{sec:reg-rom}, we first present the classical L-ROM and then introduce the novel ADL-ROM.
In Section~\ref{sec:numerical-results}, we perform a numerical investigation of the new ADL-ROM.
Specifically, we compare the ADL-ROM with the classical L-ROM and the standard G-ROM in two numerical tests: the Burgers equation with a small diffusion coefficient and the convection-dominated flow past a backward-facing step. 
Finally, in Section~\ref{sec:conclusions}, we present our conclusions and outline future research directions.




\section{G-ROM}
    \label{sec:g-rom}


In this section, we briefly describe the construction of the standard G-ROM.    
As a mathematical model, we consider the incompressible NSE: 
\begin{eqnarray}
\frac{\partial \bu}{\partial t} - Re^{-1} \Delta\bu +\bigl(\bu\cdot\nabla\bigr)\bu+\nabla p \,&=\, \bff, & \text{in } \Omega \times (0, T]\label{eqn:nse-1} \\
\nabla\cdot\bu\,&=\,0, & \text{in } \Omega \times (0, T], \label{eqn:nse-2}
\end{eqnarray}
where 
$\Omega$ is the spatial domain, $T$ is the final time, $\bu = [u_1,u_2,u_3]^\top: \Omega \times [0, T] \to \mathbb{R}^3$ is the velocity vector field, $p: \Omega \times [0, T] \to \mathbb{R}$ the pressure field, $Re$ the Reynolds number, and $\bff: \Omega \times [0, T] \to \mathbb{R}^3$ the forcing term. 
The NSE are equipped with appropriate boundary and initial conditions. 

To build the G-ROM, we first collect the snapshots $\{ \bu_h^0, \bu_h^1, \ldots, \bu_h^M\}$, which are FOM solutions at the time instances $t_0 \equiv 0, t_1, \ldots, t_M \equiv T$, which are assumed to be equispaced for simplicity of presentation.
Further details about the FOM are postponed to Section~\ref{sec:numerical-results}, as they are not essential for the presentation of the G-ROM.
Next, we compute the centered snapshots $\{ \bu_h^0 - \bU_h, \bu_h^1 - \bU_h, \ldots, \bu_h^M - \bU_h\}$, where the centering trajectory $\bU_h$ of the flow is defined as
$\bU_h(\bx) = \frac{1}{M + 1} \sum_{k=0}^M \bu_h^k(\bx)$. 
Finally, we use the proper orthogonal decomposition (POD)~\cite{HLB96,volkwein2013proper} on the centered snapshots 
to construct the ROM basis functions $\bphi_1, \ldots, \bphi_r$, where $r$ is the ROM dimension.
Although in our numerical experiments we use the POD to construct the ROM basis, we note that other ROM bases could be used instead~\cite{brunton2019data,hesthaven2015certified,quarteroni2015reduced,RozzaStabileBallarin2022}.
In what follows, we assume that the ROM velocity approximation can be written as 
\begin{eqnarray}
\bu_r(\bx,t) \, 
=\,\bU_h(\bx)+\sum_{j=1}^rc_j(t)\bphi_j(\bx)\, ,\label{eqn:rom-soln}
\end{eqnarray}
where 
$\bc(t) = [c_1(t),\cdots,c_r(t)]^\top$ are the sought ROM coefficients.
The next step in the G-ROM construction is to replace $\bu$ with $\bu_r$ in~\eqref{eqn:nse-1} and project the resulting equations onto the space spanned by the ROM basis, $\{\bphi_j\}_{j=1}^r$.
This yields the G-ROM:
\begin{eqnarray}
\left(\frac{\partial \bu_r}{\partial t},\bphi_i \right)
+\bigl((\bu_r\cdot\nabla)\bu_r,\bphi_i\bigr)
+ Re^{-1} \bigl(\nabla\bu_r,\nabla\bphi_i\bigr) = \bigl(\bff,\bphi_i\bigr)\,,
\ \ i =1,\cdots,r\, , \label{eqn:g-rom}
\end{eqnarray}
where $( \cdot , \cdot )$ denotes the $L^2$ inner product.
We notice that the pressure does not appear anymore in \eqref{eqn:g-rom} since the ROM basis functions $\bphi_1, \ldots, \bphi_r$ are divergence free.  Alternative ROMs which also allow to recover the pressure are discussed, e.g., in \cite{ballarin2015supremizer,kean2020error,noack2005need,rozza2007stability,stabile2018finite}. Since the Reg-ROMs introduced in Section~\ref{sec:reg-rom} only differ in the treatment of the nonlinear term (which does not involve the pressure), dropping the pressure unknown from the ROM is beneficial in view of a simpler exposition of the novel Reg-ROMs.
The G-ROM can be written as the following 
system of differential equations for the vector of time coefficients, $\bc = \bc(t)$:
\begin{eqnarray}
	\bc' 
	= \bb+A  \bc 
	+ \bc^\top B  \bc \, , 
	\label{eqn:g-rom-U}
\end{eqnarray}
where the vector $\bb$, the matrix $\bA$, and the tensor $\bB$ are defined as follows: 
\begin{align}
&\bb_i \,=\, \bigl(\bphi_i,\bff\bigr)-\bigl(\bphi_i,\bU_h\cdot\nabla\bU_h\bigr)
- Re^{-1} \bigl(\nabla\bphi_i,\nabla\bU_h\bigr)\, ,\label{eq:vectorb}\\[0.3cm]
&\bA_{im} \,=\, -\bigl(\bphi_i,\bU_h\cdot\nabla\bphi_m\bigr) - \bigl(\bphi_i,\bphi_m\cdot\nabla\bU_h\bigr)
- Re^{-1} \bigl(\nabla\bphi_i,\nabla\bphi_m\bigr)\, ,\\[0.3cm]
&\bB_{imn} \, =\, -\bigl(\bphi_i,\bphi_m\cdot\nabla\bphi_n\bigr)\, ,
\end{align}
for $i, n, m = 1, \hdots, r$.

As illustrated in the numerical investigation in Section~\ref{sec:numerical-results}, although G-ROM~\eqref{eqn:g-rom} works well for laminar flows, for under-resolved convection-dominated flows it yields numerical oscillations.

\section{Approximate Deconvolution (AD) ROM}
    \label{sec:ad-rom}

In this section, we introduce the AD strategy, which is a new concept in ROM stabilization.
AD is central in image processing and inverse problems~\cite{bertero1998introduction,hansen2010discrete}.
In a nutshell, the AD goal can be formulated as follows:
Given an approximation of the filtered input signal, $\obu \doteq G \bu$, where $G$ is an invertible spatial filter, find an approximation of the unfiltered input signal, $\bu$.
Of course, the exact deconvolution, $G^{-1} \obu$, would seem a natural choice.
Computationally, however, the exact deconvolution is a very bad idea since the noise in the high wavenumber components of $\obu$ is amplified by the inverse filter, $G^{-1}$.
Thus, in practice, AD strategies are used instead~\cite{bertero1998introduction,hansen2010discrete}. 
In large eddy simulations of turbulent flows, the AD models have been pioneered by Adams and Stolz for classical numerical discretizations~\cite{SA99}.
In reduced order modeling, the 
only AD model was proposed in~\cite{xie2017approximate}, where AD was used to develop a ROM closure model.

In this paper, we leverage the AD strategy to develop a different type of ROM: the ADL-ROM (Section~\ref{sec:adl-rom}), which is a ROM stabilization.
To construct the new ADL-ROM, we first present the ROM differential filter (Section~\ref{ROM-DF}).
Then, we present three different AD strategies: 
the van Cittert AD approach (Section~\ref{sec:van-cittert}), the Tikhonov AD approach (Section~\ref{sec:tikhonov}), and the Lavrentiev AD approach (Section~\ref{sec:lavrentiev}).
We note that the Lavrentiev AD approach was used in~\cite{xie2017approximate} to develop a ROM closure model.
To our knowledge, however, the van Cittert and Tikhonov approaches have never been used to develop ROM AD operators (see~\cite{cordier2010calibration,wang20162d,weller2009robust} for Tikhonov methods for ROM regularizations).
We also note that, in Section~\ref{sec:bfstep}, we perform a preliminary numerical investigation of the three AD strategies.

In what follows, we will use the following expansions:
\begin{align}
    \bu_{r}(\bx,t) &= \bU_h(\bx) + \sum_{i=1}^{r} c_i(t) \boldsymbol{\varphi}_i(\bx), 
    \qquad
    \bc \doteq [ c_1, \ldots, c_{r}]^{\top}, 
    \label{eqn:expansion-u} \\
    \overline{\bu}_{r}(\bx,t) &= \bU_h(\bx) + \sum_{i=1}^{r} \overline{c}_i(t) \boldsymbol{\varphi}_i(\bx), 
    \qquad
    \obc \doteq [ \oc_1, \ldots, \oc_{r}]^{\top}, 
    \label{eqn:expansion-ou} \\
    \bu_{r}^{AD}(\bx,t) &= \bU_h(\bx) + \sum_{i=1}^{r} c^{AD}_i(t) \boldsymbol{\varphi}_i(\bx),  
    \quad
    \bc_{AD} \doteq [ c^{AD}_1, \ldots, c^{AD}_{r}]^{\top} 
    \label{eqn:expansion-uad} .
\end{align}
That is, each of the different functions $\bu_r$ (unfiltered ROM velocity), $\overline{\bu}_r$ (filtered ROM velocity), and $\bu_{r}^{AD}$ (ROM approximate deconvolution) 
is a linear combination of the same ROM basis functions $\bphi_1, \ldots, \bphi_r$ with different coefficients. 

\begin{remark}
To simplify the presentation, in the rest of this section we assume that $\bU_h = \boldsymbol{0}$ and that the NSE~\eqref{eqn:nse-1}-\eqref{eqn:nse-2} are equipped with homogeneous Dirichlet boundary conditions. Thus, the centering trajectory, $\bU_h$, will drop from \eqref{eqn:expansion-u}-\eqref{eqn:expansion-uad}. We note, however, that a nonzero centering trajectory will actually be required in Section \ref{sec:bfstep} to handle the nonzero inlet conditions. The methodology introduced in this section readily extends to that case, which yields additional terms on the right-hand side of \eqref{eqn:rom-df-linear-system}, \eqref{eqn:CittertMatrices}, \eqref{eqn:TikhonovMatrices}, and \eqref{eqn:LavrentievMatrices}. 
For instance, 
following the same approach as that used in \eqref{eq:vectorb}, 
we obtain the following modification of \eqref{eqn:rom-df-linear-system}:
\begin{eqnarray*}
    \left( \bM + \delta^2 \bS \right) \obc
    = \bM \bc + \boldsymbol{g},
\end{eqnarray*}
where
\begin{align*}
&\boldsymbol{g}_i \,=\, -\bigl(\bphi_i,\bU_h\bigr)
- \delta^2 \bigl(\nabla\bphi_i,\nabla\bU_h\bigr)\,.
\end{align*}
\end{remark}

\subsection{The ROM Differential Filter}\label{ROM-DF}

To develop the ADL-ROM, we will use the \emph{ROM differential filter (ROM-DF)}: 
\begin{align}
    \overline{\bu} 
    &= G \bu 
    \doteq (I - \delta^2 \Delta)^{-1} \bu.       
    \label{eqn:Filter}
\end{align}
That is, given a function $\bu$, we filter it 
with the ROM-DF $G = (I - \delta^2 \Delta)^{-1}$ to obtain a filtered function $\overline{\bu}$. 
A filter radius, $\delta$, appears in the definition of $G$, and is set by the user. The input $\bu$ will typically be the velocity of the NSE~\eqref{eqn:nse-1}-\eqref{eqn:nse-2}.
We note that ROM-DF~\eqref{eqn:Filter} is equivalent to solving the PDE
\begin{eqnarray}
    && (I - \delta^2 \Delta) \obu
    = \bu
    \qquad \text{in } \Omega, 
    \label{eqn:rom-df-pde} \\
    && \obu= \boldsymbol{0}
    \hspace*{2.7cm} \text{on } \partial \Omega,  
    \nonumber 
\end{eqnarray}
where $\partial \Omega$ is the boundary of the spatial domain, $\Omega$.

We notice that equations \eqref{eqn:Filter}-\eqref{eqn:rom-df-pde} are reported in strong form, and for any input function $\bu$. In practice, however, 
we apply 
this filter to the ROM velocity, $\bu_r$, and do so by solving the ROM linear system corresponding to the application of the Galerkin method for \eqref{eqn:rom-df-pde} on the space spanned by the ROM basis functions $\{\bphi_1, \ldots, \bphi_r\}$: 
\begin{eqnarray}
    \left( \bM + \delta^2 \bS \right) \obc
    = \bM \bc,
    \label{eqn:rom-df-linear-system}
\end{eqnarray}
where $\bM$ is the POD mass matrix (which is the identity matrix since the POD basis functions are orthonormal) and $\bS$ is the POD stiffness matrix, 
whose entries are defined as follows: 
\begin{align*}
&\bM_{im} \,=\, \bigl(\bphi_i,\bphi_m\bigr)\, , \qquad \bS_{im} \,=\, \bigl(\nabla\bphi_i,\nabla\bphi_m\bigr)\, ,
\qquad i, m = 1, \hdots, r.
\end{align*}

The goal is now to obtain ${\bu}_{r}^{AD} \approx G^{-1} \overline{\bu}_{r}$ without resorting to the inverse operator $G^{-1}$, but applying instead an AD operator. Indeed, 
applying directly the filter inversion would result in an algorithm characterized by poor conditioning and ill-posedness 
\cite{bertero1998introduction,hansen2010discrete,xie2015ms}. 
Three 
AD operators are introduced below.

\subsection{The Van Cittert AD}
    \label{sec:van-cittert}
Given an integer $N$, which represents the order of the AD operator, the van Cittert approximate deconvolution operator is given by 
\begin{equation}
    \bu_{AD} = D_N \obu = \sum_{n=0}^N (I - G)^n \overline{\bu},\label{eqn:CitEqn}
\end{equation}
presented in \cite[section 3.1.1]{layton2012approximate}.
To solve this, we use a Richardson iteration, described in \cite[section 3.3.4]{layton2012approximate}
\begin{align}
    \bu_{AD}^{(n + 1)} = \bu_{AD}^{(n)} + \{\obu - G \bu_{AD}^{(n)}\}, \qquad n = 0, \hdots, N - 1, \label{eqn:CittertIteration}
\end{align}
where the initial $\bu_{AD}^{(0)}$ is set as $\overline{\bu}$. By running the iteration~\eqref{eqn:CittertIteration} $N$ times, we obtain  $\bu_{AD} \doteq \bu_{AD}^{(N)}$ as the $N^{\text{th}}$ order van Cittert operator applied to $\overline{\bu}$. However, to make this iteration practical, we need to recall that each application of the filter $G$ requires the inversion of a differential operator. Therefore, denoting by $\widetilde{\bu}^{(n + 1)}$ the quantity $\widetilde{\bu}^{(n + 1)} \doteq G \bu_{AD}^{(n)}$, and substituting in our filter from 
\eqref{eqn:Filter}, we 
observe that $\widetilde{\bu}^{(n + 1)}$ is found by solving
\begin{align}
    (I - \delta^2 \Delta)^{-1} \bu_{AD}^{(n)} = \widetilde{\bu}^{(n + 1)} \Longleftrightarrow
    (I - \delta^2 \Delta)  \widetilde{\bu}^{(n + 1)} = \bu_{AD}^{(n)}. \label{eqn:CittertEqn}
\end{align}
Equation \eqref{eqn:CittertEqn} requires to solve a linear system at each iteration $n = 0, \hdots, N - 1$. As discussed above, in order to employ the van Cittert AD in a ROM setting, by multiplying~\eqref{eqn:CittertEqn} by each test function in our ROM space and expanding our prospective solution as a linear combination of ROM basis functions, 
we obtain the linear system 
\begin{equation}
    \left( \bM + \delta^2 \bS \right) \widetilde{\bc}^{(n+1)} = \bM \bc_{AD}^{n}.\label{eqn:CittertMatrices}
\end{equation}
Thus, the iterative process~\eqref{eqn:CittertIteration} 
amounts to setting $\bc_{AD}^{(0)} = \obc$, updating the coefficients of the ROM AD velocity as
\begin{align*}
    \bc_{AD}^{(n + 1)} = \bc_{AD}^{(n)} + \{\obc - \widetilde{\bc}^{(n+1)}\} \qquad n = 0, \hdots, N - 1, 
\end{align*}
and finally defining $\bc_{AD}\doteq\bc_{AD}^{(N)}$.

\subsection{The Tikhonov AD}
    \label{sec:tikhonov}
The Tikhonov method of approximate deconvolution is defined as
\begin{equation}
\bu_{AD} = D_{\mu}^{T} \obu = (G^*G + \mu I)^{-1} G^* \overline{\bu},\label{eqn:TikEqn}
\end{equation}
where $G^*$ denotes the adjoint of the operator $G$, and $\mu \in \mathbb{R}^+$ is a positive constant; we refer, e.g., to \cite[section 3.3.1]{layton2012approximate} for further details on the Tikhonov AD. When plugging in the specific filter from 
\eqref{eqn:Filter} and proceeding formally, we can write
\begin{align}
    [G^*G + \mu I] \bu_{AD} &= G^* \overline{\bu}, \nonumber\\
    [ (I - \delta^2 \Delta)^{-*} (I - \delta^2 \Delta)^{-1} + \mu I] \bu_{AD} &= (I - \delta^2 \Delta)^{-*} \overline{\bu} ], \nonumber\\
  [(I - \delta^2 \Delta)^{-1}  + \mu (I - \delta^2 \Delta)^{*}]  \bu_{AD} &= \overline{\bu}, \nonumber\\
   [I + \mu (I - \delta^2 \Delta)(I - \delta^2 \Delta)^{*}] \bu_{AD} &= (I - \delta^2 \Delta) \overline{\bu}, \nonumber\\
    [I + \mu (I - \delta^2 \Delta^* - \delta^2 \Delta + \delta^2 \Delta \delta^2 \Delta^*)] \bu_{AD} &= (I - \delta^2 \Delta) \overline{\bu}, \nonumber \\
    \left[ I + \mu(I - 2\delta^2 \Delta + \delta^4  \Delta \Delta ) \right] \bu_{AD} &= (I - \delta^2 \Delta)\overline{\bu}. \label{eq:TikExplain}
\end{align}
In~\eqref{eq:TikExplain}, we have used the fact that the Laplacian is self-adjoint. 
In the numerical comparison of the three AD approaches in Section~\ref{sec:bfstep}, for ease of implementation in our FE setting, we neglect the term $\delta^4  \Delta^2 $.
Thus, expanding our desired solution as linear combinations of ROM basis functions, we solve for the vector of coefficients $\bc_{AD}$ satisfying  
\begin{equation}
    [(1+\mu) \bM + 2 \mu \delta^2 \bS] \bc_{AD} = (\bM + \delta^2 \bS) \overline{\bc}. \label{eqn:TikhonovMatrices}
\end{equation}

\subsection{The Lavrentiev AD}
    \label{sec:lavrentiev}
The Lavrentiev method~\cite{xie2017approximate}, which is the AD strategy we use in our numerical investigation in Section~\ref{sec:numerical-results}, is a modification of the Tikhonov method suitable for self-adjoint, positive definite operators. The filter described in \eqref{eqn:Filter} has these properties, as shown in \cite{xie2017approximate}. Given $\mu \in \mathbb{R}^+$, the Lavrentiev AD operator is given in \cite[section 3.3.2]{layton2012approximate} as  
\begin{equation}
    {\bu}_{AD} = D_{\mu}^{L} \obu = (G + \mu I)^{-1} \overline{\bu}. \label{eqn:LavrentEqn}
\end{equation}
We solve for ${\bu}_{AD}$ by substituting in $G$ 
\begin{align}
    {\bu}_{AD} &= \left[ (I - \delta^2 \Delta)^{-1} + \mu I \right]^{-1} \overline{\bu}, \nonumber \\
    \left[ (I - \delta^2 \Delta)^{-1} + \mu I \right] {\bu}_{AD} &= \overline{\bu}, \nonumber \\
    \left[I + \mu (I - \delta^2 \Delta)\right] {\bu}_{AD} &= (I - \delta^2 \Delta) \overline{\bu}, \nonumber \\
    \left[ (1 + \mu)I - \mu \delta^2 \Delta \right] {\bu}_{AD} &= (I - \delta^2 \Delta) \overline{\bu}. \nonumber
\end{align}
Expanding our solution as a linear combination of ROM basis functions, multiplying by a test function, and integrating by parts, we obtain the system  
\begin{equation}
    \left[ (1+\mu) \bM + \mu\delta^2 \bS\right] \bc_{AD} = \left(\bM + \delta^2 \bS\right) \overline{\bc} \label{eqn:LavrentievMatrices}
\end{equation}
and solve for the vector of coefficients $\bc_{AD}$.

\section{Reg-ROMs}
    \label{sec:reg-rom}

In this section, we present the Reg-ROMs that will be used in the numerical investigation in Section~\ref{sec:numerical-results}.
First, in Section~\ref{sec:l-rom}, we present the standard L-ROM.
Then, in Section~\ref{sec:adl-rom}, we introduce the novel ADL-ROM.

\subsection{L-ROM}
    \label{sec:l-rom}

The initial use of the Leray model in 1934 by Jean Leray \cite{leray1934sur} was aimed at providing a theoretical framework for establishing the existence of weak solutions for the NSE. More recently, researchers have employed the Leray model as a computational tool for simulating convection-dominated flows (e.g., turbulent flows), utilizing conventional numerical methods~\cite{geurts2003regularization}, such as the finite element method~\cite{layton2012approximate}.
We note that, when the differential filter is used, the Leray model is similar to the NS-$\alpha$ model advocated by Foias, Holm, Titi, and their collaborators~\cite{foias2001naviersteokesalpha}. 
During the last decade, the Leray model was extended to the realm of reduced order modeling, first for the Kuramoto-Sivashinsky  equations~\cite{sabetghadam2012alpha}, and then for the NSE~\cite{wells2015phd,wells2017evolve} and the quasi-geostrophic equations~\cite{Girfoglio2023CRM,Girfoglio2023JCAM}.

In our setting, the L-ROM time discretization reads as follows: Given $\boldsymbol{u}_{r}^{n}$ and $\boldsymbol{u}_{r}^{n-1}$, find $\boldsymbol{u}_{r}^{n+1}$ such that 

\begin{equation}\label{LROMeq}
	\left(\frac{\boldsymbol{u}_{r}^{n+1}-\frac{4}{3}\boldsymbol{u}_{r}^{n}+\frac{1}{3}\boldsymbol{u}_{r}^{n-1}}{\Delta t}, \boldsymbol{\varphi}_{i}\right)+\frac{2}{3} Re^{-1} \left(\nabla\boldsymbol{u}_{r}^{n+1}, \nabla \boldsymbol{\varphi}_{i}\right)+\frac{2}{3}\left(\left(\overline{\boldsymbol{u}}_{r}^{n+1} \cdot \nabla\right) \boldsymbol{u}_{r}^{n+1}, \boldsymbol{\varphi}_{i}\right)=0,
\end{equation}
$\forall n=1, \ldots, M - 1$ and $\forall i=1, \ldots, r$, where $\Delta t$ is the time step. 
Utilizing the ROM-DF~\eqref{eqn:Filter} outlined in Section~\ref{ROM-DF}, the
filtered convective term in (\ref{LROMeq}) is defined as follows: 
\begin{equation}\label{LROMuf}
	\overline{\boldsymbol{u}}_{r}^{n+1}(\boldsymbol{x}, t) \doteq \sum_{j=1}^{r} \overline{c}_{j}^{n+1}(t) \boldsymbol{\varphi}_{j}(\boldsymbol{x}) .
\end{equation}
The coefficients $\overline{\bc}^{n+1}$ in (\ref{LROMuf}) are found by solving the following reduced linear system: 
\begin{equation}\label{LROMsys}
	\lp \bM + \delta^2 \bS \rp \overline{\bc}^{n+1} = \bM\bc^{n+1},
\end{equation}
where $\bM$ and $\bS$ are the POD mass and stiffness matrices, respectively, $\delta$ is the radius of the ROM-DF, and $\bc^{n+1}$ is the vector of reduced coefficients of the input variable,  $\boldsymbol{u}_r^{n+1}$.
We note that linear system~\eqref{LROMsys} is linear system \eqref{eqn:rom-df-linear-system} for input coefficients $\bc^{n+1}$.
We also note that a second-order backward differences (BDF2) time discretization was used in (\ref{LROMeq}), but other time discretizations are possible.

Thus, at the time instance $t^{n+1}$, 
the L-ROM time discretization yields
the following nonlinear system:
\begin{equation}\label{LROMeqmat}
    \frac{1}{\Delta t} \bM \left( \bc^{n+1} - \frac{4}{3} \bc^{n} + \frac{1}{3} \bc^{n-1}\right) + \frac{2}{3}  Re^{-1}  \bS \bc^{n+1} + \frac{2}{3} \boldsymbol{C}(\overline{\bc}^{n+1}) \bc^{n+1} = {\bf 0},
\end{equation}
where $\bc^{n+1}$ 
is the unknown reduced coefficient vector of the unfiltered solution at 
time $t^{n+1}$, 
$\overline{\bc}^{n+1}$ the filtered velocity coefficient vector at time $t^{n+1}$, 
and the matrix 
$\boldsymbol{C}(\overline{\bc}^{n+1})$ is defined as
\begin{equation}
    \boldsymbol{C}(\overline{\bc}^{n+1})_{ij} = \left( (\overline{\bu}^{n+1}_r \cdot \nabla)\boldsymbol{\varphi}_i, \boldsymbol{\varphi}_j \right).
\end{equation}
System \eqref{LROMeqmat} is a nonlinear system of equations in the unknowns $\bc^{n+1}$, which we solve by 
using the Newton method. In particular, at each iteration of the Newton method, given the current approximation to $\bc^{n+1}$, we update the vector $\overline{\bc}^{n+1}$ via \eqref{LROMsys}, and evaluate the residual of \eqref{LROMeqmat}. We notice that the additional linear system \eqref{LROMsys} to be solved at each nonlinear iteration is still of small dimension $r \times r$. 

In Algorithm~\ref{alg:LROM}, we outline the main steps of the L-ROM discretization.
\begin{algorithm}
	\caption{L-ROM Pseudocode}\label{alg:LROM}
	\begin{algorithmic}[1]
        \State $\boldsymbol{u}_{-1}$, $\boldsymbol{u}_0$, $\boldsymbol{u}_{in}$, $r$   \Comment{Inputs needed}
        \For{$n \in \{1, \dots, M - 1\}$} \Comment{Time loop}
            \State \text{FOM simulation to compute $\bu_h^{n+1}$}\Comment{Snapshot collection}
        \EndFor
        \State $\mathbf{U}^r \doteq \text{POD}\left( \left\{\boldsymbol{u}_h^n\right\}_{n=1}^{M}; r\right)$ \Comment{
        POD for velocity space}
        \For{$n \in \{1, \dots, M - 1\}$} \Comment{Time loop}
            \State \text{Solve} (\ref{LROMeqmat}) to compute $\bc^{n+1}$ which requires to \Comment{L-ROM}
            \State \hspace{4.0pt} \text{solve} (\ref{LROMsys}) at each nonlinear iteration \Comment{Leray filtering}
        \EndFor
	\end{algorithmic}
\end{algorithm}

L-ROM is a significant improvement over the standard G-ROM~\eqref{eqn:g-rom} since spatial filtering alleviates the spurious numerical oscillations displayed by G-ROM in convection-dominated flows. 
However, when the ROM filter radius, $\delta$, is too large (i.e., the ROM filtering is too aggressive), L-ROM can be overdiffusive and can yield inaccurate results, as shown in the numerical investigation in Section~\ref{sec:numerical-results}.
Furthermore, numerical investigations have revealed that L-ROM is sensitive with respect to $\delta$: small $\delta$ variations can significantly decrease the L-ROM accuracy.

\subsection{ADL-ROM}
    \label{sec:adl-rom}

To address these L-ROM drawbacks, in this section we propose a novel Reg-ROM, the {\it approximate deconvolution Leray ROM (ADL-ROM)}, which is inspired from earlier work in the finite element setting based on a clever idea of Adrian Dunca~\cite{layton2008numerical,layton2008high,layton2012approximate}.
The new ADL-ROM is based on a simple yet powerful idea: 
Replace the filtered L-ROM velocity, $\obu$, with an approximate deconvolution velocity, $D(\obu)$, where $D$ is one of the AD operators defined in Section~\ref{sec:ad-rom}.
The numerical investigation in Section~\ref{sec:numerical-results} shows that using the AD operator in the ADL-ROM increases the accuracy of the L-ROM nonlinear term without compromising the L-ROM stability.
Thus, in a sense, ADL-ROM can be regarded as a  compromise between the numerically unstable G-ROM and the potentially overdiffusive L-ROM.

In our setting, the ADL-ROM time discretization reads as follows: Given $\boldsymbol{u}_{r}^{n}$ and $\boldsymbol{u}_{r}^{n-1}$, find $\boldsymbol{u}_{r}^{n+1}$ such that 

\begin{equation}\label{ADLROMeq}
	\left(\frac{\boldsymbol{u}_{r}^{n+1}-\frac{4}{3}\boldsymbol{u}_{r}^{n}+\frac{1}{3}\boldsymbol{u}_{r}^{n-1}}{\Delta t}, \boldsymbol{\varphi}_{i}\right)+\frac{2}{3} {Re^{-1}}  \left(\nabla\boldsymbol{u}_{r}^{n+1}, \nabla \boldsymbol{\varphi}_{i}\right)+\frac{2}{3}\left(\left({\mathrm{D}_{\mu}^{L}}(\overline{\boldsymbol{u}}_{r}^{n+1}) \cdot \nabla\right) \boldsymbol{u}_{r}^{n+1}, \boldsymbol{\varphi}_{i}\right)=0,
\end{equation}
$\forall n=1, \ldots, M-1$ and $\forall i=1, \ldots, r$.

Utilizing the ROM differential filter~\eqref{eqn:Filter} outlined in Section~\ref{ROM-DF} and the Lavrentiev AD operator~\eqref{eqn:LavrentEqn} in Section~\ref{sec:ad-rom}, the
AD convective term in (\ref{ADLROMeq}) is defined as follows:
\begin{equation}\label{ADLROMuf}
	{\boldsymbol{u}}_{AD, \ r}^{n+1}(\boldsymbol{x}, t) \doteq {\mathrm{D}_{\mu}^{L}}(\overline{\boldsymbol{u}}_{r}^{n+1}) \doteq \sum_{j=1}^{r} {c}_{AD, \ j}^{n+1}(t) \boldsymbol{\varphi}_{j}(\boldsymbol{x}).
\end{equation}
The coefficients ${\bc}_{AD}^{n+1}$ in (\ref{ADLROMuf}) are found by solving the reduced linear system in (\ref{eqn:LavrentievMatrices}), upon providing the current approximation of the filtered coefficients ${\obc}^{n+1}$, which in turn are computed by solving the linear system \eqref{LROMsys}.
Although we utilize the BDF2 method for the ADL-ROM time discretization~\eqref{ADLROMeq}, other time discretizations are possible.

Thus, at the time instance $t^{n+1}$, 
the ADL-ROM time discretization yields the following nonlinear system:
\begin{equation}\label{ADLROMeqmat}
    \frac{1}{\Delta t} \bM \left( \bc^{n+1} - \frac{4}{3} \bc^{n} + \frac{1}{3} \bc^{n-1}\right) + \frac{2}{3}  {Re^{-1}}  \bS \bc^{n+1} + \frac{2}{3} \boldsymbol{C}({\bc}_{AD}^{n+1}) \bc^{n+1} = {\bf 0},
\end{equation}
where $\bc^{n+1}$ are the unknown reduced coefficient vectors of the unfiltered  velocity field, and ${\bc}^{n+1}_{AD}$ the vectors of the corresponding AD ROM velocity field 
obtained by solving \eqref{LROMsys} and \eqref{eqn:LavrentievMatrices}. Therefore, each iteration of the Newton solver requires to solve two small $r \times r$ linear systems before the evaluation of the residual. 

In Algorithm~\ref{alg:ADLROM}, we outline the main steps of the ADL-ROM discretization.
\begin{algorithm}
	\caption{ADL-ROM Pseudocode}\label{alg:ADLROM}
	\begin{algorithmic}[1]
        \State $\boldsymbol{u}_{-1}$, $\boldsymbol{u}_0$, $\boldsymbol{u}_{in}$, $r$   \Comment{Inputs needed}
        \For{$n \in \{1, \dots, M - 1\}$} \Comment{Time loop}
            \State \text{FOM simulation to compute $\bu_h^{n+1}$} \Comment{Snapshot collection}
        \EndFor
        \State $\mathbf{U}^r \doteq \text{POD}\left( \left\{\boldsymbol{u}_h^n\right\}_{n=1}^{M}; r\right)$ \Comment{
        POD for velocity space}
        \For{$n \in \{1, \dots, M - 1\}$} \Comment{Time loop}
            \State \text{Solve} (\ref{ADLROMeqmat}) to compute $\bc^{n+1}$ which requires to \Comment{L-ROM}
            \State \hspace{4.0pt} \text{solve} (\ref{LROMsys}) at each nonlinear iteration, and to\Comment{Leray filtering}
            \State \hspace{4.0pt} \text{solve} (\ref{eqn:LavrentievMatrices}) at each nonlinear iteration\Comment{AD}
            
        \EndFor
	\end{algorithmic}
\end{algorithm}

\section{Numerical Results}
To maintain consistency with the only previous time AD methods were used in ROMs \cite{xie2017approximate} and for reasons discussed in Section \ref{sec:bfstep}, the rest of this paper will focus on the Lavrentiev AD method described in Section \ref{sec:lavrentiev}. Further references to the ADL-ROM in our numerical results should be understood to refer to ADL-ROMs constructed with the Lavrentiev regularization method. 
    \label{sec:numerical-results} 

\subsection{Burgers Equation}
\paragraph{Goals}
The following is an investigation of ADL-ROM performance in the simpler setting of Burgers equation, which is a common test problem for new ROM methods, see, e.g.,~\cite{ahmed2018stabilized,koc2023bounds,KV99,KV01}.

\paragraph{Computational Setting}
Burgers equation is a nonlinear problem of one spatial dimension  given as  
\begin{equation}
    u_t - \nu u_{xx} + u u_x = 0. \label{eqn:Burgers}
\end{equation}
Here, $u \doteq u(x,t)$ is defined for $x \in \Omega = [0,1]$, $t \in (0,1]$, and we consider the low viscosity setting where $\nu$ is small. A typical form of Burgers equation to make it more challenging in a computational setting is to use a discontinuous initial condition: 
\begin{equation}
    u(x,0) = \begin{cases}
    1, & 0 < x < 0.5, \\ \label{eqn:Burgers_IC}
    0, & x = 0, \ 0.5 \leq x \leq 1.
\end{cases}
\end{equation}

Because of the similarity in structure between Burgers equation and the NSE, it may be solved at a FE and ROM level using methods similar to those discussed in Section \ref{sec:g-rom}, upon neglecting the pressure term in \eqref{eqn:nse-1} and the divergence constraint \eqref{eqn:nse-2}.

\paragraph{Criteria} We will consider the errors between a well resolved FE solution $u(t)$ and a ROM solution $u_r(t)$ as 
\begin{equation}
    \label{eqn:burgerserror}
        E_{u}^{a}(t) \doteq \left\|u(t)-u_{r}(t)\right\|_{L^{2}(\Omega)}.
\end{equation}
As a comparison between the ADL-ROM and the standard G-ROM, we will also employ a relative reduction in error defined as 
\begin{equation}
 RE(t) =   -100 \cdot \frac{\left(E_{u}^{a}(t)_{G-ROM} - E_{u}^{a}(t)_{ADL-ROM}\right)}{ E_{u}^{a}(t)_{G-ROM}}. \label{eqn:burgersRELerror}
\end{equation}
 With this formula, negative percent values indicate that the ADL-ROM is performing better than the G-ROM.

\begin{figure}[H]
  \centering
  \subfloat[FEM Solution]{\includegraphics[width=0.49\textwidth]{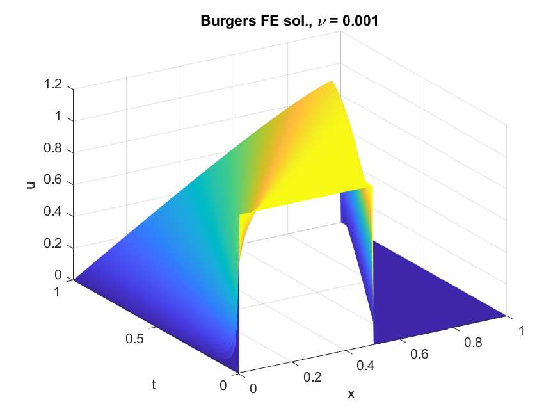}\label{fig:FEM_FULL}}
  \hfill
  \subfloat[G-ROM, 10 basis functions]{\includegraphics[width=0.49\textwidth]{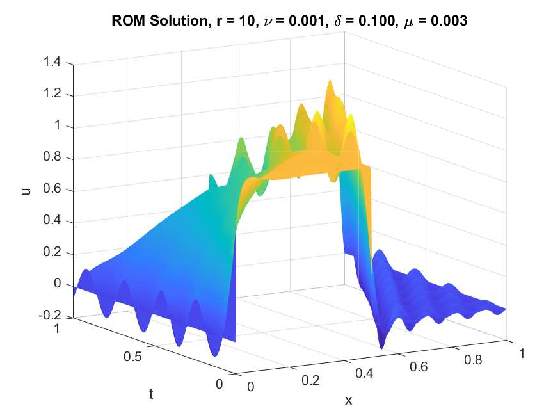}\label{fig:G-ROM_FULL}}
  \caption{ Burgers equation~\cite{moore2023filter-ms}:
  (a) FEM solution;
  (b) G-ROM solution.
  }
\end{figure}

\paragraph{Numerical results}
The plots in Fig.~\ref{fig:FEM_FULL}--\ref{fig:G-ROM_FULL} show that, even with a well resolved FE solution, the nonlinear nature of the problem as well as the challenging initial condition present problems for the standard G-ROM. 
The FE solution in Figure~\ref{fig:FEM_FULL} was  computed with $350$ evenly spaced elements between $0$ and $1$ and piecewise linear finite element functions. Both the FE solution in Fig.~\ref{fig:FEM_FULL} and the G-ROM solution in Fig.~\ref{fig:G-ROM_FULL} were computed with an implicit Euler time discretization and using a Newton solver, 75 even time steps between 0 and 1, and $\nu = 10^{-3}$. Even though the FE solution is completely resolved, the ROM displays severe oscillations. As shown below, using instead an ADL-ROM constructed with the Lavrentiev regularization method, we are able to lessen the oscillations. 

\par 

\begin{figure}[H]
  \centering
  \subfloat[ROM solutions for $\delta = 0.5$.]{\includegraphics[width=0.49\textwidth]{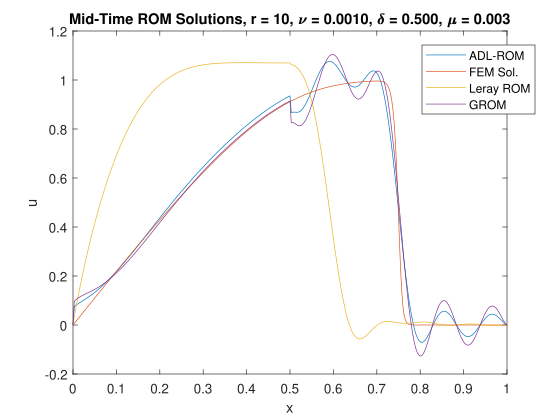}\label{fig:delta05compare}}
  \hfill
  \subfloat[ROM solutions for $\delta = 0.1$.]{\includegraphics[width=0.49\textwidth]{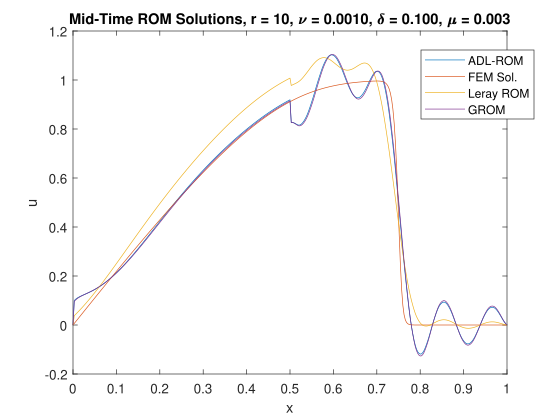}
   \label{fig:delta01compare}}
  \caption{Comparison of ROM solutions at $t = 0.5$ with two $\delta$ values.}
  \label{fig:ADL_WELL_Resolved}
\end{figure}

We emphasize that, even though we are using a smoothing process, it is important to avoid over smoothing. The objective is to minimize non-physical phenomena, not replace one (oscillations) with another (excessive smoothing). In Figure \ref{fig:delta05compare}, the Leray ROM with a large filter radius of $\delta = 0.5$ produces a very smooth solution which is incompatible with the FE solution, while the ADL-ROM tapers the oscillations present in the G-ROM. In Figure \ref{fig:delta01compare}, a more moderate value of $\delta = 0.1$ produces a L-ROM which is still over-smoothed, though the ADL-ROM is now nearly identical to the G-ROM. 
\par 

The plots in Figure \ref{fig:ADL_WELL_Resolved} demonstrate that the ADL-ROM is able to accomodate a larger filter radius than the L-ROM. Our investigation suggests that the larger values of $\delta$, near $\delta = 0.5$, are more suitable for minimizing our definition of error given in~\eqref{eqn:burgerserror}. Additionally, the dramatic differences between the L-ROM in Figure \ref{fig:delta05compare} and Figure \ref{fig:delta01compare} as compared to the ADL-ROM solutions show that the addition of the parameter $\mu$ causes the ADL-ROM solutions to vary less dramatically with changes in $\delta$.

\begin{figure}[H]

  \centering
  \subfloat[$L^2$ errors per time step]{\includegraphics[width=0.5\textwidth]{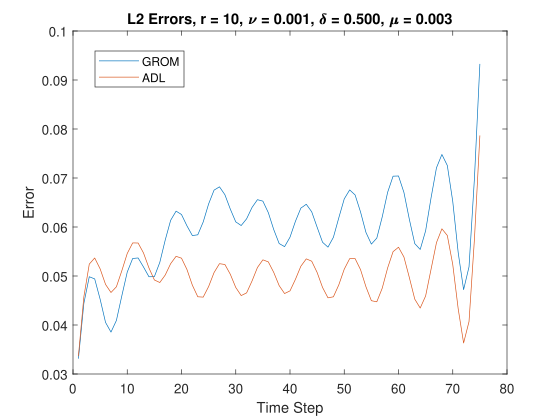}\label{fig:ADL_WELL_ERRORS_L2}}
  \hfill
  \subfloat[Relative \% error change (Eqn. \eqref{eqn:burgersRELerror})]{\includegraphics[width=0.5\textwidth]{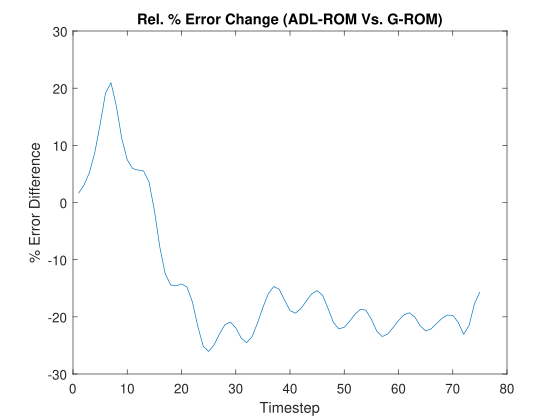}\label{fig:ADL_WELL_ERRORS_PERC}}
  \caption{ADL-ROM errors \cite{moore2023filter-ms}: $\delta = 0.5$, $\mu = 0.003$. }
  \label{fig:ADL_Error_WELL}
\end{figure}

While the effects of the ADL-ROM shown in Figure \ref{fig:delta05compare} are mild at each time step, the errors in Figure \ref{fig:ADL_WELL_ERRORS_L2} and Figure \ref{fig:ADL_WELL_ERRORS_PERC} show that a reasonable reduction in error can be obtained from the ADL-ROM, which is only a single additional linear solve from the L-ROM. In Table \ref{tab:Burgers_AVG_l2_error_FIRST}, the mean error of the L-ROM are of order $10^{-1}$, while the other two methods are of order $10^{-2}$, so the L-ROM errors for each time step are not plotted in Figure \ref{fig:ADL_WELL_ERRORS_L2} to avoid distorting the scale.
\begin{table}[ht]
    \centering
    \begin{tabular}{|c|c|c|c|} \hline
     Model &  G-ROM  & L-ROM & ADL-ROM \\ \hline
    Mean Error &  $5.94 \times 10^{-2}$ & $4.24 \times 10^{-1}$ &  $5.03 \times 10^{-2}$ \\
    \hline
    \end{tabular}
    \caption{Mean absolute $L^2$ errors across all time steps.}
    \label{tab:Burgers_AVG_l2_error_FIRST}
\end{table}
These tests should not be taken to suggest that the L-ROM cannot be tuned provide acceptable results, but that the ADL-ROM is a minimal cost adaption to the L-ROM, the main benefit of which is the ability to choose a larger radius filter than is possible with the L-ROM. The second parameter $\mu$ is an additional complexity in choice, but at the same time the solutions generated by the ADL-ROM are less variable than those generated by the L-ROM. 

\begin{figure}[H]
  \centering
  \subfloat[FE solution]{\includegraphics[width=0.49\textwidth]{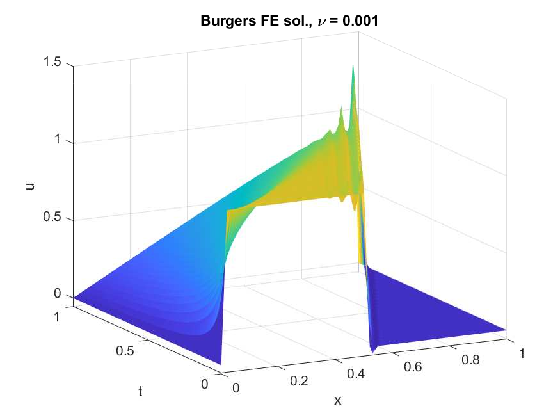}\label{fig:FEM_50}}
  \hfill
  \subfloat[G-ROM, 10 basis functions]{\includegraphics[width=0.49\textwidth]{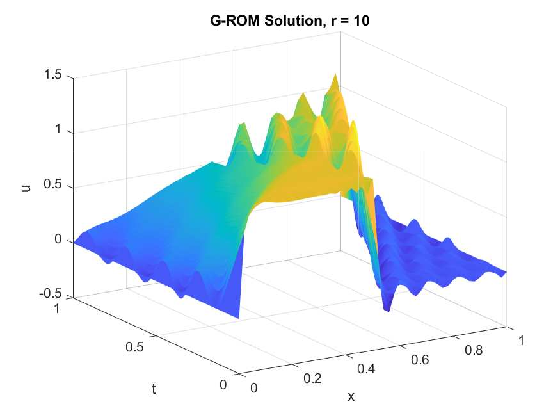}\label{fig:G-ROM_50}}
  \caption{The G-ROM propagates the oscillations present in the FE solution.}
  \label{fig:poorly_resolved_burgers}
\end{figure}

\paragraph{Additional Test}
Following on the observation that the ADL-ROM allows larger filter values of $\delta$, the following test aims to increase the numerical oscillations in the solution and assess the ADL-ROM's performance in a more challenging setting. If we take the same initial condition as in~\eqref{eqn:Burgers_IC} but decrease the number of evenly spaced elements from 350 to 50, the shock will fail to be properly resolved in the FE simulation and cause cascading oscillations in Figure \ref{fig:FEM_50}, which are propagated into the G-ROM of Figure \ref{fig:G-ROM_50}. The images in Figure \ref{fig:poorly_resolved_burgers} were generated with $\nu = 10^{-3}$, 100 evenly spaced timesteps between 0 and 1, but only 50 evenly spaced elements. 

\begin{figure}[H]
  \centering
  \subfloat[ADL-ROM solution]{\includegraphics[width=0.49\textwidth]{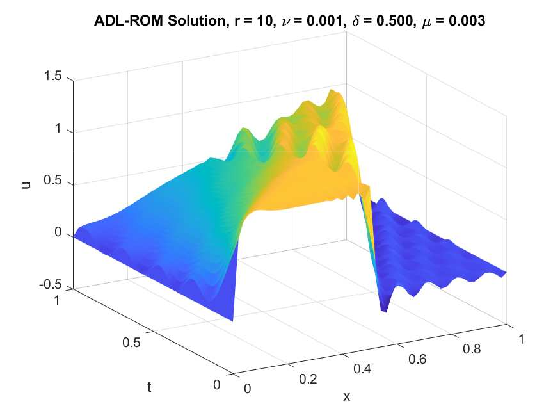}\label{fig:ADL_50}}
  \hfill
  \subfloat[Relative \% error change (Eqn. \eqref{eqn:burgersRELerror})]{\includegraphics[width=0.49\textwidth]{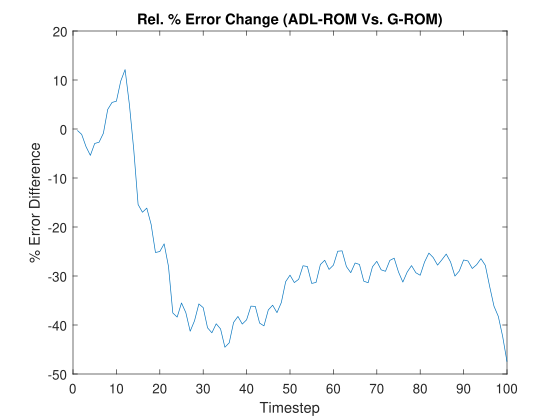}\label{fig:50_rez_percent}}
  \caption{The ADL-ROM damps the oscillations.}
  \label{fig:burgers_50_rez_ADL}
\end{figure}

Because the input FE data is poor, this necessitates a change in how the errors in Figure \ref{fig:50_rez_percent} are calculated. Instead of comparing against the poor input FE data, we compute a higher resolution FE solution.
In Figure \ref{fig:50_rez_percent}, the reduction in error is computed by using a FE solution with 200 evenly spaced elements and 100 timesteps. 
\noindent

With ADL parameters chosen as $\delta = 0.5$, $\mu = 0.003$, the ADL-ROM shown in Figure \ref{fig:burgers_50_rez_ADL} is able to substantially quell the oscillations inherent in the FE and G-ROM solutions. Comparing between Tables $\ref{tab:Burgers_AVG_l2_error_FIRST}$ and $\ref{tab:Burgers_AVG_l2_error}$, the ADL-ROM is able to maintain its accuracy significantly better than the G-ROM is when downgrading the input data. The leading digits of error for the L-ROM are the same in both cases because the L-ROM fails in the same way in both cases by severely over smoothing the computed solution. 

\begin{table}[ht]
    \centering
    \begin{tabular}{|c|c|c|c|} \hline
     Model &  G-ROM  & L-ROM & ADL-ROM \\ \hline
    Mean Error &  $8.05 \times 10^{-2}$ & $4.24 \times 10^{-1}$ &  $5.68 \times 10^{-2}$ \\
    \hline
    \end{tabular}
    \caption{Mean absolute $L^2$ errors across all time steps.}
    \label{tab:Burgers_AVG_l2_error}
\end{table}
\vspace{2mm}

These results suggest that the ADL-ROM performs best when it is applied as a treatment for significant oscillations in the numerical solution. In particular, this suggests the applicability of the method for use in convection-dominated flow simulations, such as the ones we propose for the NSE.

\subsection{Backward Facing Step}\label{sec:bfstep}

    \paragraph{Goals} In this section, we analyze and compare the performances of G-ROM (\ref{eqn:g-rom}), L-ROM (\ref{LROMeq}), and ADL-ROM (\ref{ADLROMeq}) for the NSE (\ref{eqn:nse-1})-(\ref{eqn:nse-2}). The goal is to investigate ADL-ROM accuracy. We consider the two-dimensional velocity flow over the backward-facing step at Reynolds number ${Re} = 1429$.

    \begin{figure}[h!]
        \centering
        \includegraphics[width=0.9\textwidth]{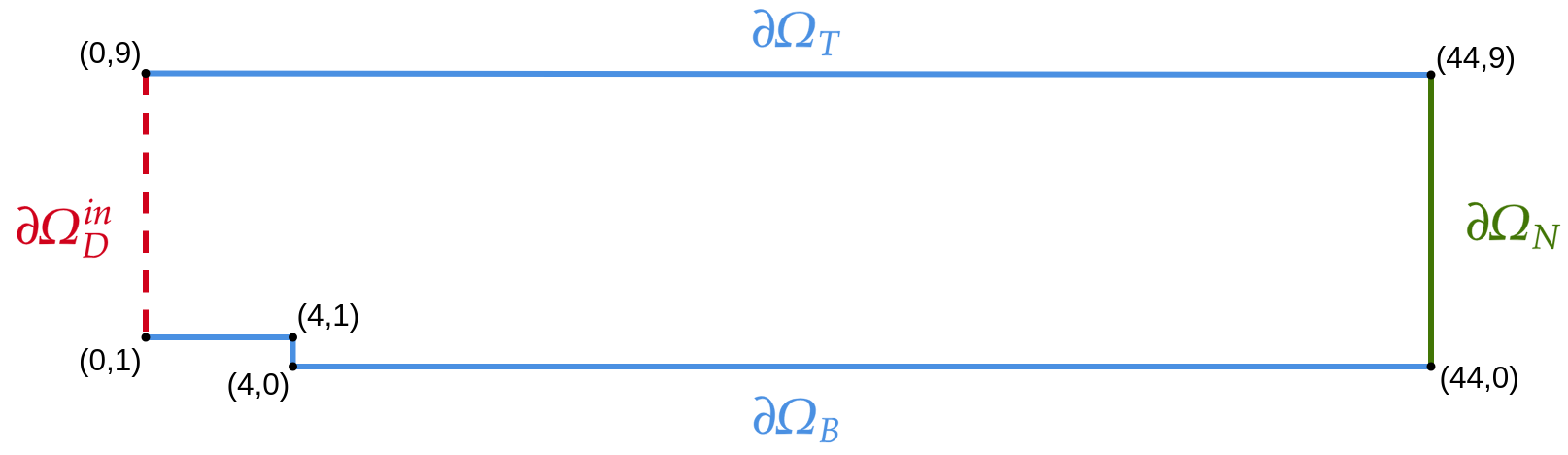}
            \caption{The computational domain, $\Omega$. $\partial \Omega_{D} \doteq \partial \Omega_{D}^{i n} \cup \partial \Omega_{D}^{\text {wall }}$, where the inlet boundary $\partial \Omega_{D}^{i n}$ is represented by a dashed red line and the no-slip boundaries by a solid blue line.}
        \label{BDF_domain}
    \end{figure}

    \paragraph{Computational setting} We consider the motion of an incompressible flow in the same domain as that used in \cite[section 4.4]{baiges2015reduced} and in \cite[section 3.5]{mou2021data}, i.e, $\Omega \doteq\{[0,44] \times[0,9]\} \backslash \{[0,4] \times[0,1]\}$, which is depicted in Figure \ref{BDF_domain}. We set a kinematic viscosity of $7\times10^{-4}$ and use no-slip boundary conditions on $\partial \Omega_{D}^{\text {wall }} \doteq \partial \Omega_{\text {B}} \cup \partial \Omega_{\text {T}}$, representing the union of the bottom $\partial \Omega_{\text {B}}$ and top walls $\partial \Omega_{\text {T}}$ of the channel (solid blue boundary in Figure \ref{BDF_domain}), with a constant inlet velocity profile $\boldsymbol{u}_{i n}=
    (1,0)$ on $\partial \Omega_{D}^{i n}$ (red dashed line in Figure \ref{BDF_domain}). Furthermore, we employ homogeneous Neumann boundary conditions on $\partial \Omega_{N}$ (green line in Figure \ref{BDF_domain}).

    \begin{figure}[h!]
        \centering
        \includegraphics[width=0.9\textwidth]{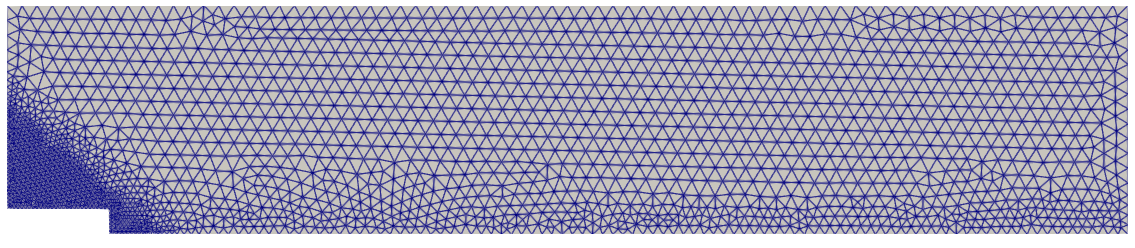}
            \caption{The FE mesh.}
        \label{BDF_mesh}
    \end{figure}

    \paragraph{Snapshot generation}  We perform our tests on a triangular mesh with $h_{\min }=0.15$ and $h_{\max }=0.6$ (Figure \ref{BDF_mesh}). For the spatial discretization, we employ the inf-sup stable Taylor-Hood $\mathbb{P}^{2}-\mathbb{P}^{1}$ FE pair for velocity and pressure, respectively, and this leads to a FE space of dimension $N_{h} \doteq N_{h}^{u}+N_{h}^{p}=18404 + 2370 =20774$. A second-order 
    BDF2 scheme is used with the time step $\Delta t = 0.05$ for both FOM and ROM time discretizations. The time interval on which FOM is performed is $[T_0, T] = [0, 150]$. The value of the initial conditions $\boldsymbol{u}_{-1}$ and $\boldsymbol{u}_{0}$ is $({\bf 0},{\bf 0})$.

    At the FOM level, since we work in the under-resolved regime, we apply a regularization strategy. Specifically, we employ the evolve-filter-relax (EFR) strategy (for details, see \cite{strazzullo2021consistency}). The EFR stabilization strategy allows us to obtain accurate results in the convection-dominated regime on the coarse mesh in Figure \ref{BDF_mesh}. 
    
    We use EFR as a stabilization strategy for the following reasons:  EFR has been widely used at the FOM level, is easy to implement, is effective, and has been used successfully in our numerical investigations~\cite{gunzburger2019evolve,moore2023filter-ms,mou2023energy,strazzullo2021consistency,wells2015phd,wells2017evolve}.  
    We emphasize, however, that other regularization strategies (e.g., Leray or even ADL, which would ensure FOM-ROM consistency~\cite{strazzullo2021consistency}) or stabilization approaches could be investigated.  
    This could be subject of future work.
    
\begin{figure}[h!]
        \centering
        \includegraphics[width=0.8\textwidth]{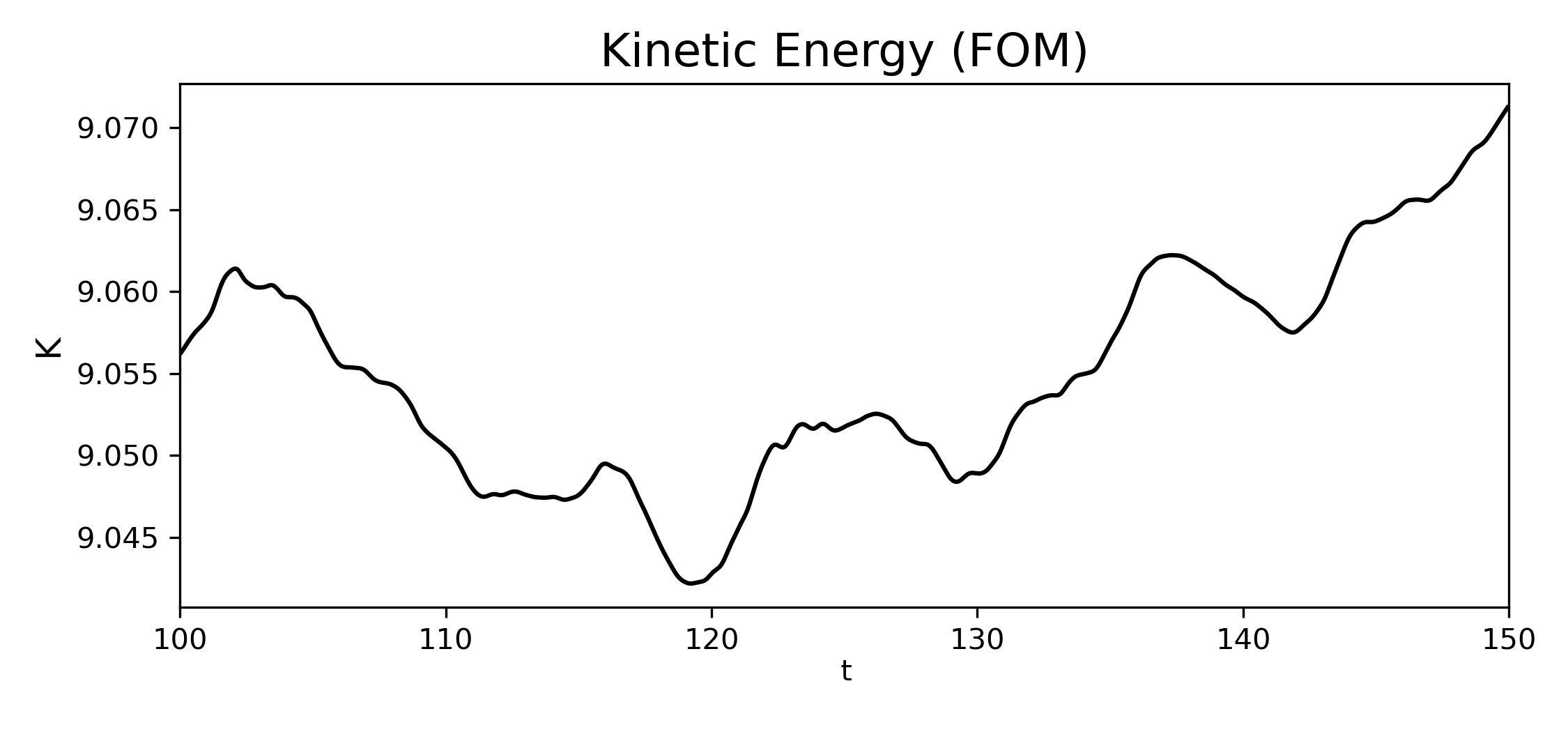}
            \caption{Time evolution of FOM kinetic energy.}
        \label{BDF_FOM_KE}
    \end{figure}

    In Figure \ref{BDF_FOM_KE}, we plot the time evolution of the FOM kinetic energy on the time interval $[100, 150]$. This plot shows that the flow 
    is not periodic or periodic-like.

    \paragraph{ROM construction} To build the ROM basis functions, we collect 1000 snapshots of the velocity field over the time interval $[100.05,150]$, as in \cite{mou2021data} and 
    \cite{reyes2020projection}. The first $r=10$ POD modes capture about $80 \% $ of the flow's energy. 
    We limit ourselves to a small number of POD basis functions for the purpose of showing the effect of regularization when the ROM cannot accurately describe all scales of motion, which is the case in many realistic applications. These modes are used 
    to construct all ROMs. All the ROMs are investigated on the 
    time interval $[100.05,150]$.

    \paragraph{AD Comparison}
    Before proceeding to the discussion of ADL-ROM results, we report here on the choice of AD to use given the options discussed in Section \ref{sec:ad-rom}. To do this, we will briefly examine the effects that these AD methods have on a standard G-ROM solution. The G-ROM solution will be calculated following the methods discussed in Section \ref{sec:g-rom} and in the previous paragraph to a final time of 20 time units with a fixed step length of 0.01 time units.  

    \begin{figure}[h!]
        \centering
        \includegraphics[width=0.9\textwidth]{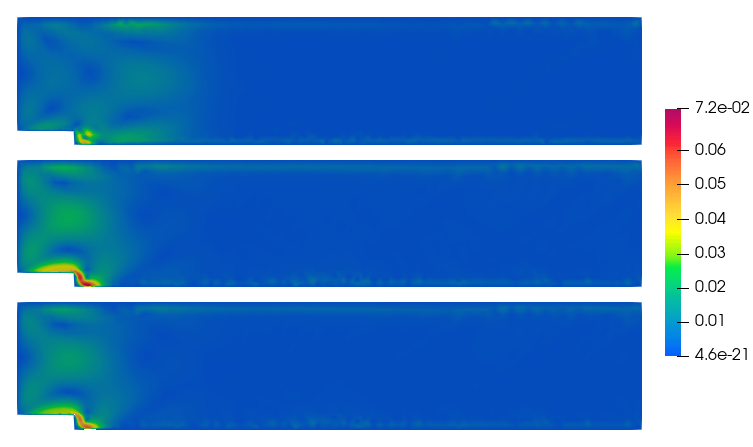}
        \caption{Pointwise difference between AD applied to G-ROM and G-ROM, $\delta = 0.6$. $1^{st}$: Lavrentiev with $\mu = 0.1$. $2^{nd}$: Tikhonov with $\mu = 0.1$. $3^{rd}$: van Cittert with $ N = 10$.}
        \label{fig:AD_Compare_Different}
    \end{figure}

\begin{figure}[h!]
        \centering
        \includegraphics[width=0.9\textwidth]{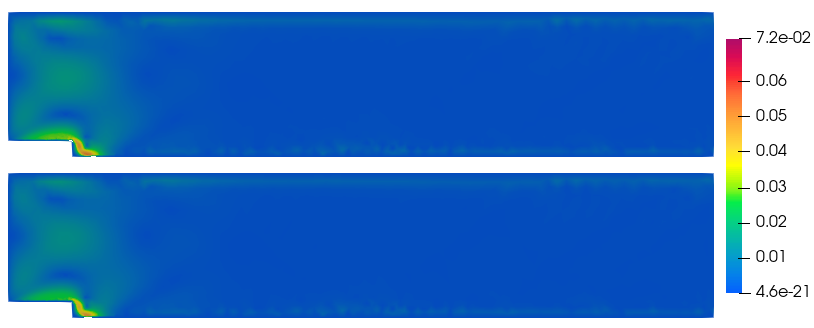}
        \caption{Pointwise difference between AD applied to G-ROM and G-ROM, $\delta = 0.6$. $1^{st}$: Tikhonov with $\mu = 0.05$. $2^{nd}$: van Cittert with $ N = 12$.}
        \label{fig:AD_Compare_Similar}
    \end{figure}

    In Figures \ref{fig:AD_Compare_Different} and \ref{fig:AD_Compare_Similar}, we plot the pointwise difference between the computed G-ROM solution at $t = 20$ and each of the AD methods applied to the G-ROM solution. Figure \ref{fig:AD_Compare_Different} also shows that, with AD parameters chosen as labeled, the Lavrentiev AD method produces results which are the least similar to the Tikhonov or van Cittert methods. However, this is mostly a matter of tuning: Figure \ref{fig:AD_Compare_Similar} shows that altering the AD parameters of $\mu$ for the Tikhonov or $N$ for the van Cittert method will make those methods produce very similar results as to the Lavrentiev method in the first plot of Figure \ref{fig:AD_Compare_Different}.
    
    Each of the AD methods shown in Figures \ref{fig:AD_Compare_Different} and \ref{fig:AD_Compare_Similar} affects the computed G-ROM solution most strongly near the backward facing step, which is the primary place where vortices are formed and is where we want the smoothing to be focused. Additionally, by varying either $\delta$ in the original filter or one of the AD parameters $\mu$ or $N$ as in Figure \ref{fig:AD_Compare_Similar}, each of the AD regularization methods can produce similar results to the other methods. For these reasons, 
    the numerical results of this paper focus on the Lavrentiev AD method described in Section \ref{sec:lavrentiev}. 

    \paragraph{Parameters} As mentioned before, the inlet velocity and the size of the step lead to a Reynolds number ${Re} = 1429$. For the EFR FOM stabilization strategy, there are two parameters: the filtering radius of the differential filter, $\delta_{EFR}$, employed in the \textit{Filter} step, and the relaxation parameter, $\chi \in [0,1]$, employed in the \textit{Relax} step. For the rest of the article, we present results for $\delta_{EFR}=0.01$ and $\chi=0.01$. This choice limits the amount of dissipation introduced by the differential filter in the EFR algorithm.

    One of the most important parameters of L-ROM and ADL-ROM is the radius of the ROM differential filter introduced in Section \ref{ROM-DF}, $\delta$.  To avoid confusion between the filtering radius of the EFR step $\delta_{EFR}$, and the filtering radius of the ROM-DF, in the discussion we will denote by $\delta_{ROM}$ the radius of the ROM-DF, instead of $\delta$. 
    To illustrate the case of an aggressive choice of the L-ROM filter radius, in our numerical investigation we used $\delta_{ROM}=0.6$.
    Another important ADL-ROM parameter 
    is the Lavrentiev regularization parameter, $\mu$. In the image processing and inverse problem communities (see \cite{bertero1998introduction,hansen2010discrete}), numerous approaches are proposed to determine the regularization parameters. In our numerical investigation, for the fixed $\delta_{ROM}=0.6$ value, we choose the $\mu$ value that ensures that the ADL-ROM solution is as close as possible to the FOM solution. This approach yields the 
    value $\mu=0.06$.

    \paragraph{Criteria} We test the ROMs accuracy by using the velocity field $L^{2}$ absolute and relative errors, and the kinetic energy. The errors are defined as
    \begin{equation}\label{BDF_L2_ERR_EQ}
        E_{\boldsymbol{u}}^{a}(t) \doteq \left\|\boldsymbol{u}(t)-\boldsymbol{u}_{r}(t)\right\|_{L^{2}(\Omega)} \quad \text { and } \quad E_{\boldsymbol{u}}^{r}(t) \doteq \frac{E_{\boldsymbol{u}}^{a}(t)}{\|\boldsymbol{u}(t)\|_{L^{2}(\Omega)}},
    \end{equation}
    while the kinetic energy is defined as
    \begin{equation}\label{BDF_KE_EQ}
        K_{\boldsymbol{u}} \doteq \frac{1}{2} \|\boldsymbol{u}(t)\|_{L^{2}(\Omega)}^{2}.
    \end{equation}
    To compare the ROMs’ performance, we also use the following error metric
     \begin{equation}\label{BDF_AV_ERR_EQ}
            \text{time-average} \ \ L^2 \ \ \text{norm}: \ \ \mathcal{E}^{a/r}(L^2) \doteq \frac{1}{M} \sum_{j=1}^M E_{\boldsymbol{u}}^{a/r}(t_j),
        \end{equation}
        where the $a$ and $r$ superscripts refer to the absolute and relative errors, respectively.

    \paragraph{Numerical results}

    \begin{table}[h!]
        \centering
        \renewcommand\arraystretch{1.2} 
    	\begin{tabular}{ |c|c|c| }
             \hline
             \multicolumn{3}{|c|}{Average $L^2$ errors} \\
             \hline
               Model & $\mathcal{E}^{a}(L^2)$ & $\mathcal{E}^{r}(L^2)$ \\
             \hline
              {G-ROM}   & 1.2668 & 7.25e-02 \\
              {L-ROM}   & 1.6061 & 8.88e-02 \\
              {ADL-ROM} & 1.0144 & 5.60e-02 \\
             \hline
        \end{tabular}
        \caption{Average absolute and relative $L^2$ errors for G-ROM, L-ROM, and ADL-ROM.}
        \label{BDF_ROM_AV_ERR}
    \end{table}
    
    \begin{figure}[h!]
        \begin{minipage}[c]{0.5\textwidth}
            \includegraphics[width=\textwidth]{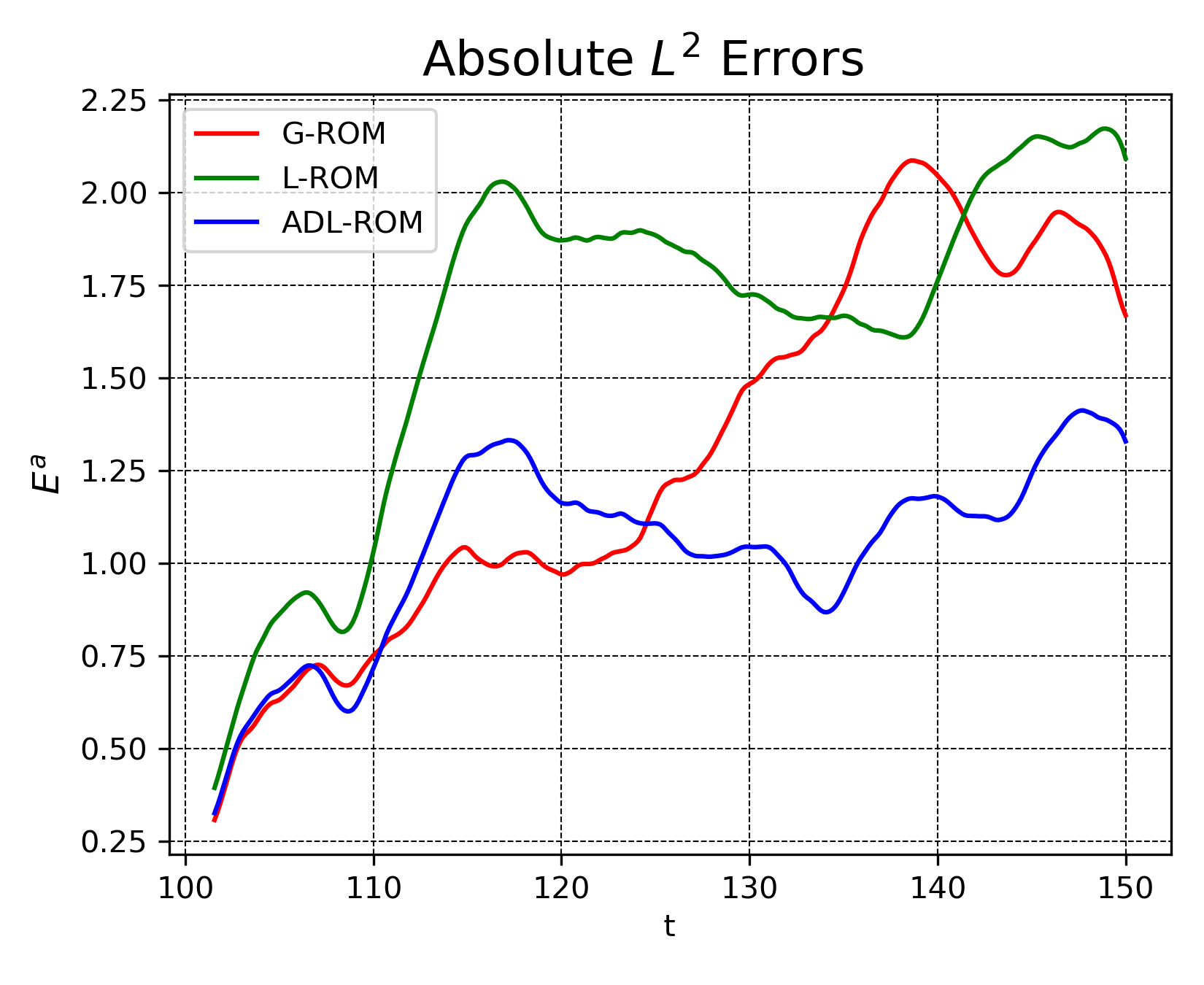}
        \end{minipage} 
        \begin{minipage}[c]{0.5\textwidth}
            \includegraphics[width=\textwidth]{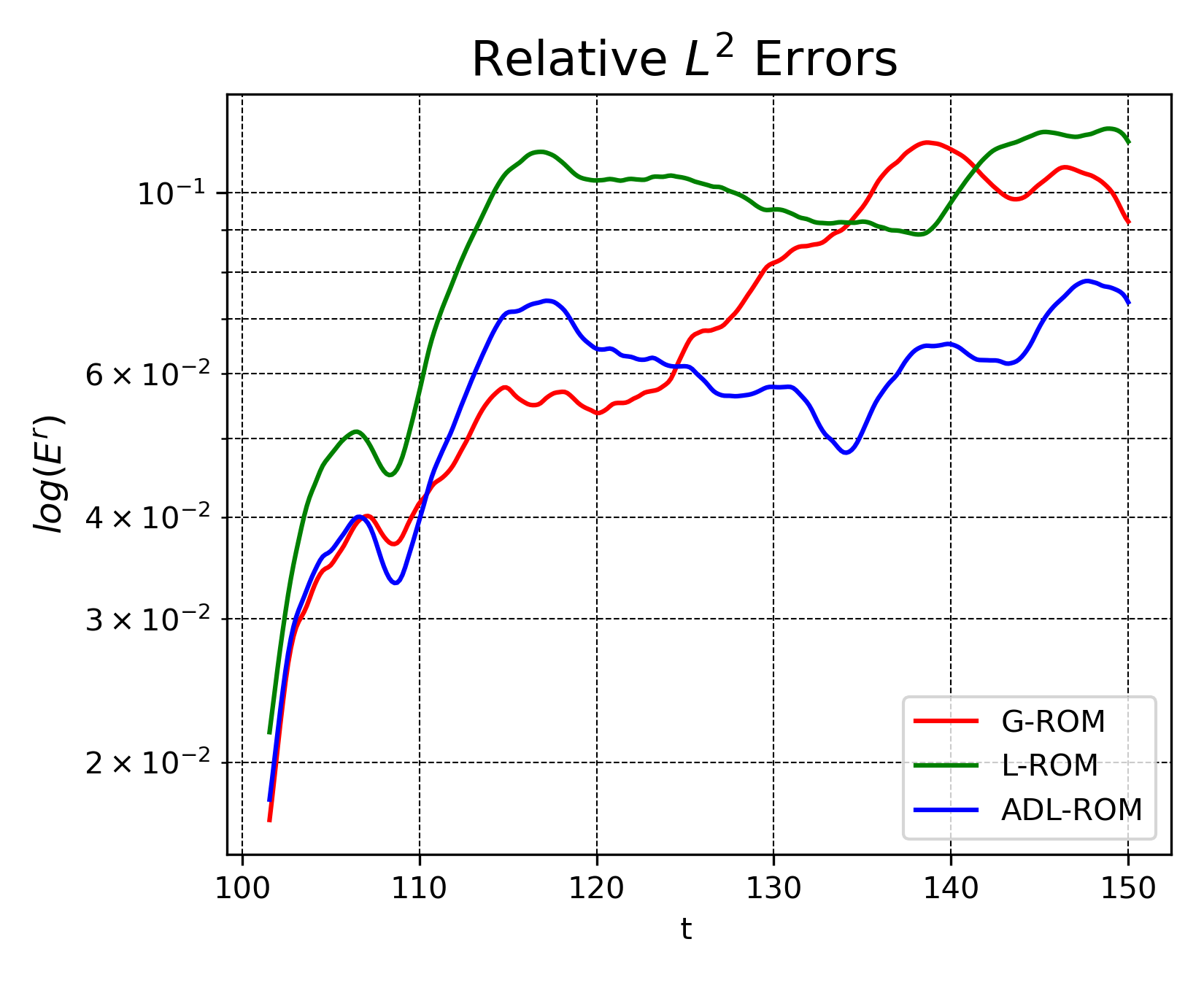}
        \end{minipage}
        \captionof{figure}{Time evolution of the absolute and relative $L^2$ errors for G-ROM (red), L-ROM (green), and ADL-ROM (blue).}
        \label{BDF_ROM_10_ERRs}
    \end{figure}

    \begin{figure}[h!]
        \centering
        \includegraphics[width=0.9\textwidth]{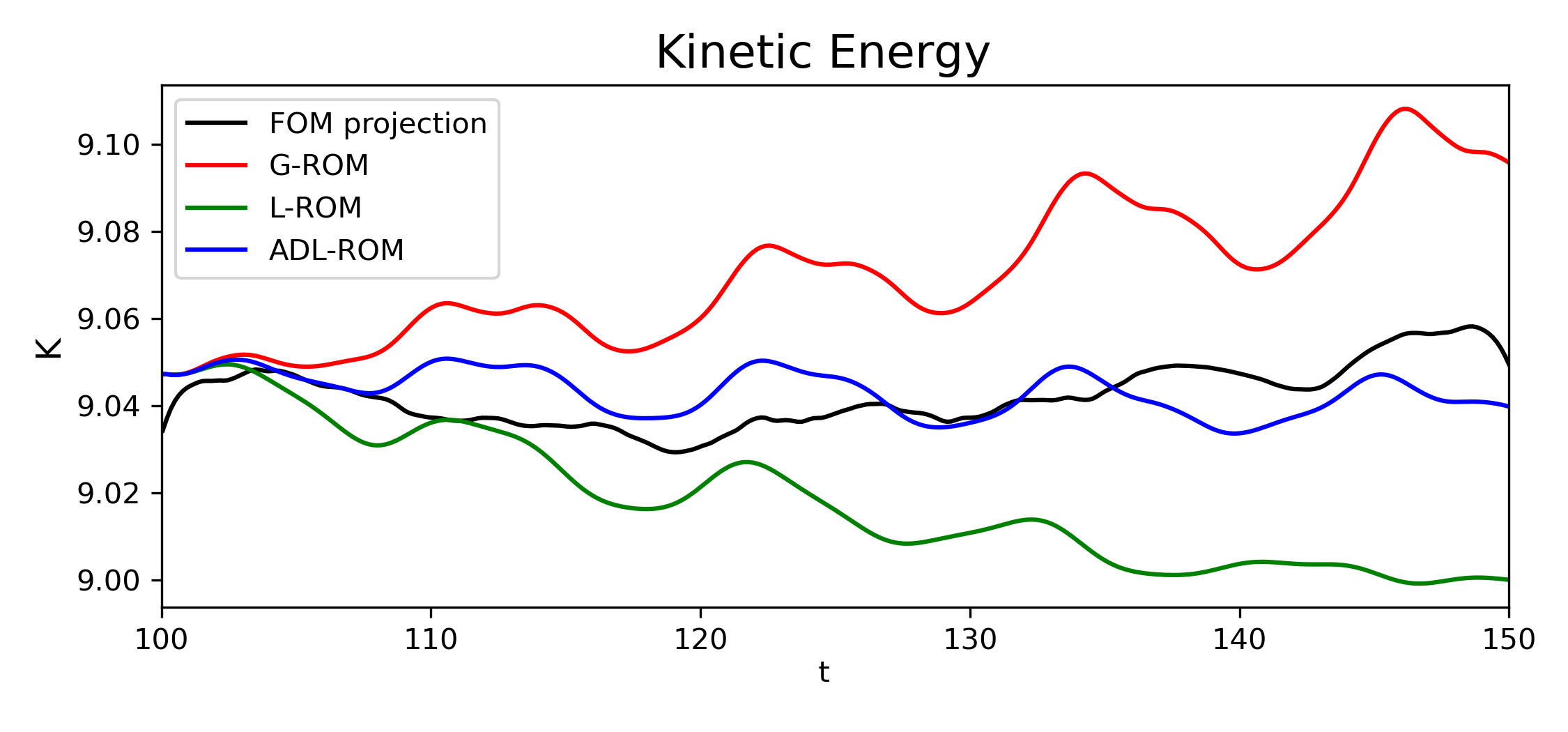}
            \caption{Time evolution of the kinetic energy for FOM projection (black), G-ROM (red), L-ROM (green), and ADL-ROM (blue).}
        \label{BDF_ROM_KE}
    \end{figure}
    
    In Table \ref{BDF_ROM_AV_ERR}, we list the average absolute and relative $L^2$ errors (\ref{BDF_AV_ERR_EQ}) for the G-ROM, L-ROM, and ADL-ROM. These results show that both the G-ROM and the ADL-ROM are more accurate than the L-ROM, with the G-ROM being about $30\%$ more accurate than the L-ROM, and the ADL-ROM about $60\%$ more accurate than the L-ROM. Furthermore, the ADL-ROM is consistently more accurate than the G-ROM with a difference of about $30\%$. 

    Also the temporal trend of the absolute and relative $L^2$ errors (\ref{BDF_L2_ERR_EQ}), displayed in Figure \ref{BDF_ROM_10_ERRs}, confirms this conclusion. It shows that ADL-ROM always performs better than L-ROM and starts to perform better than G-ROM after a while. On the other hand, G-ROM performs  better than L-ROM except for the time interval $[134, 142]$.
    In Figure \ref{BDF_ROM_KE}, we plot the time evolution of the kinetic energy (\ref{BDF_KE_EQ}) of the FOM projection, G-ROM, L-ROM, and ADL-ROM. These plots support the conclusions in Table \ref{BDF_ROM_AV_ERR} and in Figure \ref{BDF_ROM_10_ERRs}. Specifically, the G-ROM and L-ROM results are relatively inaccurate, while the ADL–ROM results are significantly more accurate.
    
\begin{figure}[h!]
        \centering
        \includegraphics[width=0.9\textwidth]{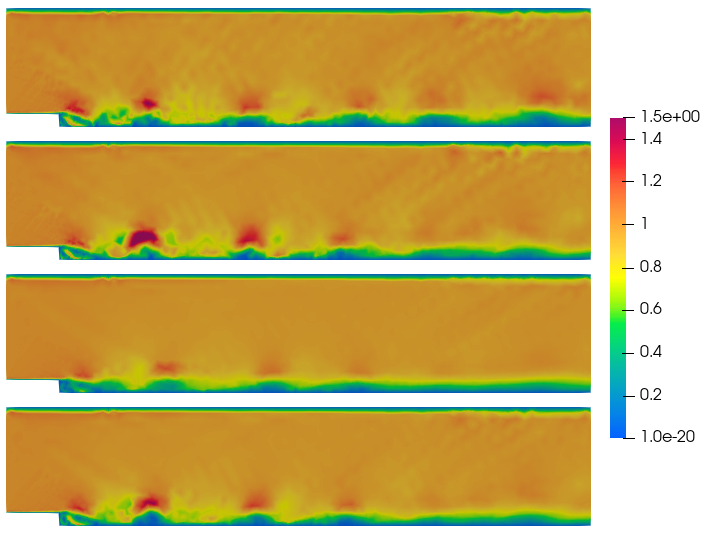}
            \caption{At $T_f=150$, $1^{st}$: FOM solution, $2^{nd}$: G-ROM solution, $3^{rd}$: L-ROM solution, and $4^{th}$: ADL-ROM solution.}
        \label{BDF_ROM_10_sol}
    \end{figure}

\begin{figure}[h!]
        \centering
        \includegraphics[width=0.9\textwidth]{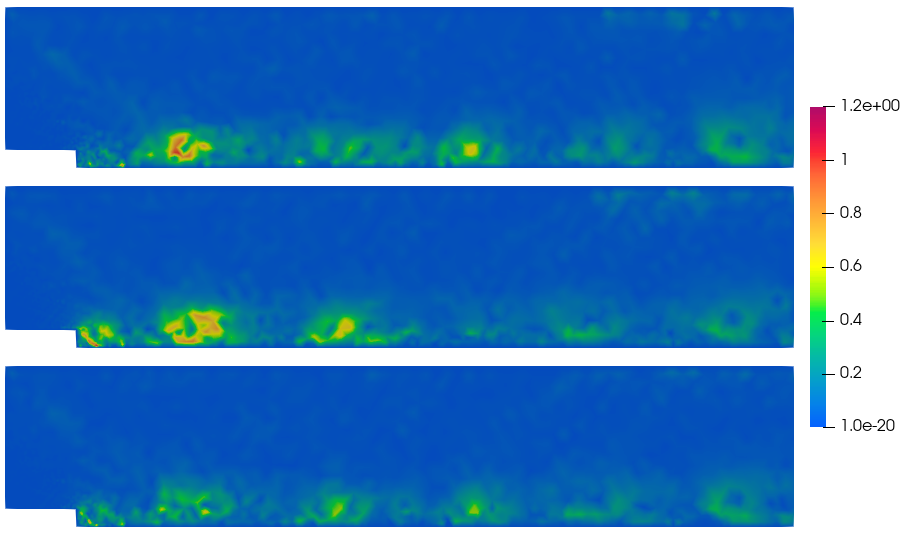}
            \caption{At $T_f=150$, $1^{st}$: pointwise error between FOM and G-ROM solutions, $2^{nd}$: pointwise error between FOM and L-ROM solutions, and $3^{rd}$: pointwise error between FOM and ADL-ROM solutions.}
        \label{BDF_ROM_10_pERR}
    \end{figure}

    We show representative solutions of the velocity field for $T_f=150$ in Figure \ref{BDF_ROM_10_sol}. It is clear that the L-ROM is not able to reconstruct the solution provided by the FOM, while the G-ROM leads to more accurate results and the ADL-ROM leads to the best results. As additional proof, we exhibit the plot of pointwise error functions, i.e., the pointwise difference between the FOM and the ROMs solutions in Figure \ref{BDF_ROM_10_pERR}. The previously observed difference in the accuracy of the ROMs solutions is even clearer in this plot by considering the maximum pointwise value in each of the three plots.

    The errors listed in Table \ref{BDF_ROM_AV_ERR} and all the plots show that both the G–ROM and the ADL–ROM are more accurate than the L-ROM. Furthermore, the ADL–ROM is more accurate than the G–ROM.

\paragraph{Computational performance}
The FE solution has $20774$ degrees of freedom. The CPU time required by the $3000$ time steps of the FE simulation is $9900 \, s$, and the CPU time required by the last $1000$ final time steps is $3000 \, s$.
The final $1000$ time steps of the ROM solutions plotted in Figure \ref{BDF_ROM_10_sol} require $2.20 \, s$ (G-ROM), $2.35 \, s$ (L-ROM), and $2.47 \, s$ (ADL-ROM).
The mild increase in CPU time from G-ROM to L-ROM is due to line 8 in Algorithm \ref{alg:LROM}, which is required only by L-ROM. Similarly, the mild increase in CPU time from L-ROM to ADL-ROM is due to line 9 in Algorithm \ref{alg:ADLROM}, which is required only by ADL-ROM. We emphasize, however, that all the ROMs yield a significant speedup: $1363.63$ (G-ROM), $1276.59$ (L-ROM) and $1214.58$ (ADL-ROM).


\section{Conclusions and Future Work}
    \label{sec:conclusions}

L-ROM is a popular ROM stabilization for convection-dominated flows~\cite{Girfoglio2019CF,girfoglio2021pod,gunzburger2020leray,iliescu2018regularized,kaneko2020towards,sabetghadam2012alpha,wells2017evolve,xie2018numerical}.
One of L-ROM's drawbacks is that its accuracy highly depends on $\delta$, the radius of the spatial filter used to smooth out the spurious numerical oscillations in the standard G-ROM.
Indeed, L-ROM numerical investigations have shown that small variations around the optimal $\delta$ values (i.e., filter radius values that yield the most accurate L-ROM results) can result in inaccurate predictions. 
For example, an aggressive $\delta$ choice, such as using a large $\delta$ value intended to alleviate as much as possible the G-ROM oscillations, generally yields inaccurate L-ROM results since in that case L-ROM is overdiffusive. 

To alleviate L-ROM's overdiffusive character in these settings, we propose a new regularized ROM, the ADL-ROM, which leverages approximate deconvolution~\cite{bertero1998introduction,layton2012approximate} to increase the L-ROM's accuracy without compromising its numerical stability.
ADL-ROM's main novelty is the replacement of L-ROM's nonlinear term $\obu \cdot \nabla \bu$ with $D(\obu) \cdot \nabla \bu$, where $D(\obu)$ is the approximate deconvolution of the filtered ROM velocity, $\obu$.
Another novel contribution of the ADL-ROM is the use of two new types of approximate deconvolution ROM operators: the van Cittert AD and the Tikhonov AD.

To assess the new ADL-ROM, we compare it with the standard L-ROM and G-ROM in the numerical investigations of two test problems: 
the Burgers equation with a small diffusion coefficient and the convection-dominated flow past a backward-facing step.
Our numerical investigation yielded the following conclusions:
First, for a large filter radius, ADL-ROM yielded more accurate results than the standard L-ROM.
Specifically, ADL-ROM added a limited amount of numerical diffusion, just enough to stabilize the ROM simulation.
L-ROM, on the other hand, was overdiffusive, yielding inaccurate results.
The second conclusion yielded by our numerical investigation was that ADL-ROM was less sensitive with respect to parameter variations than L-ROM.
Specifically, finding the optimal ADL-ROM parameters was significantly easier than finding the optimal L-ROM parameters.
The third conclusion of our numerical investigation was that, for carefully chosen parameters, the three ROM approximate deconvolution strategies investigated (i.e., Lavrentiev AD, Tikhonov AD, and van Cittert AD) yielded similar results.

This first application of approximate deconvolution to construct novel regularized ROMs such as the new ADL-ROM has yielded encouraging results.
There are, however, several research directions that could be further investigated.
For example, the novel ADL-ROM can be assessed in the numerical simulation of more challenging convection-dominated flows, such as the turbulent channel flow~\cite{mou2023energy}.
One could also further investigate the role of the AD strategy in the ADL-ROM construction, both in terms of accuracy and parameter sensitivity.
In particular, a thorough investigation of the accuracy of the new AD ROM strategies (Tikhonov and van Cittert) with ROM-DF and other filters should be performed.
Finally, providing mathematical support (such as numerical analysis and parameter scalings) for the novel ADL-ROM is an open problem.

\section*{Acknowledgments}
We acknowledge the European Union's Horizon 2020 research and innovation program under the Marie Skłodowska-Curie Actions, grant agreement 872442 (ARIA).
FB acknowledges the INdAM-GNCS project ``Metodi numerici per lo studio di strutture geometriche parametriche complesse'' (CUP E53C22001930001), the PRIN 2022 PNRR project ``ROMEU: Reduced Order Models for Environmental and Urban flows'', and the project ``Reduced order modeling for numerical simulation of partial differential equations'' funded by Università Cattolica del Sacro Cuore. TI acknowledges support through National Science Foundation grants DMS-2012253 and CDS\&E-MSS-1953113.

\bibliographystyle{plain}
\bibliography{traian}
\end{document}

%% file: notation.tex
\newcommand{\lp}{\left(}
\newcommand{\rp}{\right)}

\newtheorem{remark}{Remark}[section]

\def\PP{{{\rm l}\kern - .15em {\rm P} }}
\def\PN2{{\PP_{N}-\PP_{N-2}}}

\newcommand{\bphi}{\boldsymbol{\varphi}}

\newcommand{\bA}{\boldsymbol{A}}
\newcommand{\bb}{\boldsymbol{b}}
\newcommand{\bB}{\boldsymbol{B}}
\newcommand{\bc}{\boldsymbol{c}}

\newcommand{\bff}{\boldsymbol{f}}

\newcommand{\bM}{\boldsymbol{M}}

\newcommand{\bS}{\boldsymbol{S}}
\newcommand{\bu}{\boldsymbol{u}}
\newcommand{\bU}{\boldsymbol{U}}

\newcommand{\bx}{\boldsymbol{x}}

\newcommand{\oc}{\overline{c}}
\newcommand{\obc}{\overline{\boldsymbol c}}

\newcommand{\obu}{\overline{\boldsymbol u}}

\definecolor{vargreen}{rgb}{0.0, 0.5, 0.0}

\newcommand{\deleted}[1]{{}}
